\documentclass{amsart}
\usepackage{amssymb,latexsym}
\usepackage{graphicx}
\newtheorem{theorem}{Theorem}
\newtheorem{lemma}{Lemma}
\newtheorem{definition}{Definition}

\newtheorem{corollary}{Corollary}

\begin{document}
\title{AN ANALOGUE OF THE L\'EVY-CRAM\'ER THEOREM FOR RAYLEIGH DISTRIBUTIONS }
\author{Thu Van Nguyen}
\address{Department of Mathematics;
         International University, HCM City;
         No.6 Linh Trung ward, Thu Duc District, HCM City;
         Email: nvthu@hcmiu.edu.vn}
\date{June 30, 2009}

\begin{abstract} In the present paper we prove that every k-dimensional Cartesian product of Kingman convolutions can be embedded into a k-dimensional symmetric convolution (k=1, 2, \ldots) and obtain an analogue of the Cram\'er-L\'evy theorem for multi-dimensional Rayleigh distributions.  
    \end{abstract}\indent
\maketitle{Keywords and phrases: Cartesian products of Kingman convolutions; Rayleigh distributions; radial characteristic functions; 

AMS2000 subject classification: 60B07, 60B11, 60B15, 60K99.}
 \section{Introduction, Notations and Preliminaries}\label{S:intro}
 In probability theory and statistics, the {\bf Rayleigh distribution} is a continuous probability distribution which is widely used to model events that occur in different fields such as medicine, social and natural sciences. A multivariate Rayleigh distribution is the probability distribution of a vector of norms of random Gaussian vectors.
 The purpose of this paper, is to introduce and study the fractional indexes multivariate Rayleigh distributions via the Cartesian product of  Kingman convolutions and, in particular, to prove an analogue of the L\'evy-Cram\'er theorem for multivariate Rayleigh distributions.

  Let $\mathcal P:=\mathcal P(\mathbb R^+)$  denote the set of all probability measures (p.m.'s) on the positive half-line   $\mathbb R^+$. Put, for each continuous bounded function f on $\mathbb R^{+}$,
 \begin{multline}\label{astKi}
\int_{0}^{\infty}f(x)\mu\ast_{1,\delta}\nu(dx)=\frac{\Gamma(s+1)}{\sqrt{\pi}\Gamma(s+\frac{1}{2})}\\
\int_{0}^{\infty}\int_{0}^{\infty}\int_{-1}^{1}f((x^2+2uxy+y^2)^{1/2})(1-u^2)^{s-1/2}\mu(dx)\nu(dy)du,
\end{multline}
  where $\mu\mbox{ and }\nu\in\mathcal P\mbox{ and }\delta=2(s+1)\geq1$ (cf. Kingman \cite{Ki} and
Urbanik \cite{U1}). The convolution algebra $(\mathcal{P},\ast_{1,\delta})$ is
 the most important example of Urbanik convolution algebras (cf Urbanik \cite{U1}). In language of the
  Urbanik convolution algebras, the {\it characteristic measure}, say $\sigma_s$, of the Kingman convolution
  has the Rayleigh density
  \begin{equation}\label{Ray}
  d\sigma_s(y)= \frac{2{(s+1)^{s+1}}}{\Gamma(s+1)}y^{2s+1}\exp{(-(s+1)y^2)}dy
  \end{equation}
with the characteristic exponent $\varkappa=2$ and the kernel
$\Lambda_s$
\begin{equation}\label{eq:Lam}
 \Lambda_s(x)= \Gamma(s+1) J_{s}(x)/(1/2x)^{s},
\end{equation}
where $J_s(x)$ denotes the Bessel function of the first kind,
\begin{equation}\label{eq:Bessel}
J_s(x):= \Sigma_{k=0}^{\infty} \frac{(-1)^k
(x/2)^{\nu+2k}}{k!\Gamma(\nu+k+1)}.
\end{equation}
  It is known (cf. Kingman \cite{Ki}, Theorem 1), that the kernel $\Lambda_s$ itself is an
ordinary characteristic function (ch.f.) of a symmetric p.m., say $F_s$, defined on the
interval [-1,1]. Thus, if $\theta_s$ denotes a random variable (r.v.) with distribution
$F_s$ then for each $t\in \mathbb R^+$,
\begin{equation}\label{eq:LamThe}
\Lambda_s(t)= E\exp{(it\theta_s)}=\int_{-1}^1\cos{(tx)}dF_s(x).\end{equation}
 Suppose that $X$ is a nonnegative r.v. with distribution $\mu\in\mathcal{P}$
 and $X$ is independent of $\theta_s$.  The {\it radial characteristic function}
 (rad.ch.f.) of $\mu$, denoted by $\hat\mu(t),$ is defined by
  \begin{equation}\label{ra.ch.f.}
\hat\mu(t) = E\exp{(itX\theta_s)} = \int_0^{\infty}
\Lambda_s(tx)\mu(dx),
\end{equation}
 for every $t\in \mathbb R^{+}$.
 The characteristic measure of the Kingman convolution $\ast_{1, \delta}$, denoted by $\sigma_s$,
has the Maxwell density function
\begin{equation}\label{Maxwell density}
\frac{d\sigma_s(x)}{dx}=\frac{2(s+1)^{s+1}}{\Gamma(s+1)}x^{2s+1}exp\{-(s+1)x^2\}, \quad(0<x<\infty).
\end{equation}
and the rad.ch.f.
\begin{equation}
\hat\sigma_s(t)=exp\{-t^2/4(s+1)\}.
\end{equation}

Let $\tilde P:=\tilde{\mathcal P}(\mathbb R^+)$ denote the class of symmetric p.m.'s on $\mathbb R^+.$ Putting, for every $G\in \mathcal P$,
\begin{equation*}
 F_s(G)=\int_0^{\infty}F_{cs} G(dc),
 \end{equation*}
  we get a continuous homeomorphism from the Kingman convolution algebra $(\mathcal{P},\ast_{1,\delta})$ onto the ordinary convolution algebra $(\tilde{\mathcal P}, \ast)$ such that
\begin{eqnarray}\label{homeomorphism1}
 F_s\{G_1\ast_{1, \delta}G_2\}&=&(F_sG_1)\ast(F_sG_2) \qquad G_1, G_2\in \mathcal P\\
 F_s\sigma_s&=&N(0, 2s+1)
 \end{eqnarray}
 which shows that every Kingman convolution can be embedded into the ordinary convolution $\ast$.
    \section{Cartesian product of Kingman convolutions}
   Denote by $ \mathbb {R}^{+k}, k=1,2,...$ the k-dimensional nonnegative cone
 of $ \mathbb {R}^{k}$ and $\mathcal{P}(\mathbb {\mathbb R}^{+k})$ the class of
all p.m.'s on $\mathbb  {\mathbb R}^{+k}$ equipped with the weak convergence. In
the sequel, we will denote the multidimensional vectors and random vectors (r.vec.'s)
and their distributions by bold face letters.

For each point z of any set $A$ let $\delta_z$ denote the Dirac measure (the unit mass) at
the point z. In particular, if $\mathbf  x=(x_1, x_2,\cdots,x_k)\in
 \mathbb R^{k+}$, then
\begin{equation}\label{proddelta}
\delta_{\mathbf {x}} =
\delta_{x_1}\times\delta_{x_2}\times \ldots\times\delta_{x_k},\quad (k\; times),
\end{equation}
where the sign $"\times"$ denotes the Cartesian product of
measures.
  We put, for $\mathbf {x} = (x_1,\cdots, x_k)\mbox{ and }\mathbf {y} =
(y_1,y_2,\cdots, y_k)\in \mathbb R^{+k},$
  \begin{equation}\label{convdeltas}\delta_{\mathbf x}\bigcirc_{s, k} \delta_{\mathbf {y}} = \{\delta_{x_1}\circ _s \delta_{y_1}\} \times\{\delta_{x_2} \circ _s\delta_{y_2}\}
\times\cdots\
 \times \{\delta_{x_k} \circ_s \delta_{y_k}\},\quad (k\; times),
\end{equation}
here and somewhere below for the sake of simplicity we denote the
Kingman convolution operation $\ast_{1,\delta}, \delta=2(s+1)\ge 1$ simply by $\circ_{s}, s\ge \frac{!}{2}.$
    Since convex combinations of p.m.'s of the form
(\ref{proddelta}) are dense in $\mathcal P(\mathbb R^{+k})$ the
relation (\ref{convdeltas}) can be extended to arbitrary p.m.'s $
\mathbf{G}_1 \mbox{  and } \mathbf{G}_2\in\mathcal{P}( \mathbb R^{+k})$.
Namely, we put
\begin{equation}\label{convF}
\mathbf {G}_1 \bigcirc_{s, k} \mathbf {G}_2 = \iint\limits_{ \mathbb R^{+k}}
  \delta_{\mathbf {x}} \bigcirc_{s, k} \delta_{\mathbf {y}} {\mathbf G}_1(d\mathbf {x}) {\mathbf G}_2(d\mathbf {y})
 \end{equation} which means that for each continuous bounded function $\phi$ defined on $\mathbb R^{+k}$
 \begin{equation}\label{convof}
\int\limits_{\mathbb R^{+k}} \phi({\mathbf z}) {\mathbf G}_1 \bigcirc_{s, k} {\mathbf G}_2  (d{\mathbf z})= \iint\limits_{ \mathbb R^{+k}}\big\{\int\limits_{\mathbb R^{+k}} \phi({\mathbf z})  \delta_{{\mathbf x}} \bigcirc_{s, k} \delta_{{\mathbf y}}(d{\mathbf z})\big\}{ \mathbf G}_1(d{\mathbf x}) {\mathbf G}_2(d{\mathbf y}).
\end{equation}
 In the sequel, the binary
operation $\bigcirc_{s, k}$ will be called {\it the k-times Cartesian
product of Kingman convolutions} and the pair $(\mathcal P( \mathbb R^{+k}), \bigcirc_{s, k})$ will be called
{\it the k-dimensional Kingman convolution algebra}. It is easy to show that the
binary operation $\bigcirc_{s, k}$ is continuous in the weak topology
which together with (\ref{astKi}) and (\ref{convF}) implies the
following theorem.
 \begin{theorem}\label{Theo:Kingmanalgebra} The pair $(\mathcal P{( \mathbb R^{+k})} ,\bigcirc_{s,  k})$
 is a commutative
topological semigroup with $\delta_{\mathbf  0}$ as the unit element.
Moreover, the operation $\bigcirc_{s, k}$ is distributive w.r.t.
convex combinations of p.m.'s in $\mathcal P( \mathbb R^{+k})$.
\end{theorem}
 \  For every ${\mathbf  G}\in\mathcal P( \mathbb R^{+k})$ the
k-dimensional rad.ch.f. $\hat{{\mathbf  G}}({\mathbf  t}), {\mathbf
t}=(t_1, t_2, \cdots t_k)\in \mathbb R^{+k},$ is defined by
\begin{equation}\label{k-ra.ch.f.}
\hat{\mathbf  G}(\mathbf  {t})=\int\limits_{\mathbb R^{+k}}
 \prod_{j=1}^{k}\Lambda_s(t_jx_j){ \mathbf  G}(\mathbf {dx}),
 \end{equation}
 where $\mathbf {x}=(x_1, x_2, \cdots x_k)\in  \mathbb R^{+k}.$
 Let $\mathbf{\Theta_s} = \{\theta_{s, 1},\theta_{s, 2}, \cdots ,\theta_{s, k}\}$, where $\theta_{s, j}$ are independent r.v.'s with the same distribution
 $F_s $.
 Next, assume that $ {\mathbf X}=\{X_1, X_2,..., X_k\}$ is a k-dimensional
r.vec. with distribution $\mathbf{G}$ and $\mathbf{X}$ is independent of
$\mathbf{\Theta}_s$. We put
     \begin{equation}\label{[Theta,X]}
[{\mathbf\Theta}_s,{\mathbf X}]=\{{\theta_{s,1} X_1, \theta_{s, 2} X_2,...,\theta_{s, k}X_k}\}.
\end{equation}
 Then, the following formula is equivalent to (\ref{k-ra.ch.f.}) (cf. \cite{Ng4})
 \begin{equation}\label{multiradchf}
 \widehat{\mathbf G}({\mathbf t})=Ee^{i<{\mathbf t},[{\mathbf\Theta_s, \mathbf X}]>},\qquad {\mathbf t}\in \mathbb R^{+k}.
 \end{equation}
The Reader is referred to Corollary 2.1, Theorems 2.3 \& 2.4 \cite{Ng4} for the  principal properties of  the above rad.ch.f.
 Given $s\ge -1/2$  define a map $F_{s, k}:  \mathcal  P(\mathbb R^{+k}) \rightarrow  \mathcal P(\mathbb R^k)$ by
 \begin{equation}\label{k-map}
 F_{s, k}({\mathbf G})=\int\limits_{\mathbb R^{+k}} (T_{c_1}F_s)\times(T_{c_2}F_s)\times  \ldots\times(T_{c_k}F_s) {\mathbf G}(d{\mathbf  c}),
 \end{equation}
 here and in the sequel, for a distribution  $\mathbf G$ of a r.vec. $\mathbf X$ and a real number c we denote by $T_c{\mathbf G}$ the distribution of $c{\mathbf X}$.
  Let us denote by $\tilde{ \mathcal P}_{s, k}(\mathbb{R}^{+k})$ the sub-class of  $\mathcal P(\mathbb R^k)$ consisted of all p.m.'s defined by the right-hand side of (\ref{k-map}).
 By virtue of (\ref{k-ra.ch.f.})-(\ref{k-map}) it is easy to prove the following theorem.
 \begin{theorem}\label{symmconvo}
The set $\tilde{ \mathcal P}_{s, k}(\mathbb{R}^{+k})$ is closed w.r.t. the weak convergence and the ordinary convolution $\big.\ast$ and the following equation holds
\begin{equation}\label{Fourier=rad.ch.f.}
\hat{\mathbf G}({\mathbf t})=\mathcal F(F_{s, k}({\mathbf G}))({\mathbf t})\qquad {\mathbf t}\in {\mathbb R^{+k}}
\end{equation}
where $\mathcal F({\mathbf K})$ denotes the ordinary characteristic function (Fourier transform)
of a p.m. ${\mathbf K}$. Therefore, for any ${\mathbf G}_1\mbox{ and } {\mathbf G}_2\in \mathbb R^{+k}$
\begin{equation}\label{convolequality}
F_{s, k}({\mathbf G}_1)\big.\ast F_{s, k}({\mathbf G}_2)=F_{s, k}({\mathbf G}_1\bigcirc_{s, k}{\mathbf G}_2)
\end{equation}
and the map $F_{s, k}$ commutes with convex combinations of distributions and scale changes
$T_c, c>0$. Moreover,
\begin{equation}\label{Gaussian-Rayleigh}
F_{s, k}({\Sigma_{s, k}})=N({\mathbf 0}, 2(s+1){\mathbf I})
\end{equation}
where $\Sigma_{s, k}$ denotes the k-dimensional Rayleigh distribution and $N({\mathbf 0}, 2(s+1){\mathbf I})
$ is the symmetric normal distribution on $\mathbb R^k \mbox{ with variance operator } R= 2(s+1)
{\mathbf I}, {\mathbf I}$ being the identity operator. 
Consequently, Every Kingman convolution algebra $\big( \mathcal P(\mathbb R^{+k}), \bigcirc_{s, k}\big)$ is
embedded in the ordinary convolution algebra $\big( \mathcal P_{s, k}(\mathbb{R}^{+k}), \big.\star \big)$ and the map $F_{s, k}$ stands for a homeomorphism.
\end{theorem}
\begin{proof}

	First we prove the equation (\ref{Fourier=rad.ch.f.}) by taking the Fourier transform of the right-hand side of (\ref{k-map}). We have, for ${\mathbf t}\in \mathbb R^k,$
\begin{eqnarray}\label{Fourier-r.ch.f}
\mathcal F(F_{s, k}({\mathbf G}))({\mathbf t})&=&\notag
\int\limits_{\mathbb R^k}\Pi_{j=1}^k \cos(t_jx_j)H_{s, k}({\mathbf G})d{\mathbf x}\\
&=&\int\limits_{\mathbb R^k}\int_{\mathbb R^{+k}}
\Pi_{j=1}^k\cos(t_jx_j)(T_{c_j}F_s (d{\mathbf x}){\mathbf G}(d{\mathbf c})\\
&= &\int\limits_{\mathbb R^{+k}} \prod_{j=1}^{k}\Lambda_s(t_jc_j) {\mathbf G}(d{\mathbf c})\notag\\
&=& \hat{\mathbf G}({\mathbf t})\notag
\end{eqnarray}
which implies that  the set  set $\tilde{ \mathcal P}_{s, k}(\mathbb{R}^{+k})$ is closed w.r.t. the weak convergence and the ordinary convolution $\big.\ast$ and, moreover the equations (\ref{convolequality}) and (\ref{Gaussian-Rayleigh}) hold.
\end{proof}


   \begin{definition}\label{k-ID}
 A p.m. ${\mathbf F} \in  \mathcal P(\mathbb R^{+k})$ is called $\bigcirc_{s, k}-$infinitely divisible
 ($\bigcirc_{s, k}-$ID), if for every m=1, 2, \ldots there exists $\mathbf F_m\in  \mathbf P(\mathbb R^{+k})$ such that
 \begin{equation}\label{kID}
{ \mathbf F}={\mathbf F}_m\bigcirc_{s, k} {\mathbf F}_m\bigcirc_{s, k}\ldots \bigcirc_{s, k}{\mathbf F}_m\quad (m\;times).
 \end{equation}
 \end{definition}
 \begin{definition}\label{stability}
 $\mathbf F$ is called stable, if for any positive
 numbers a and b there exists a positive number c such that
 \begin{equation}\label{k-stability}
  T_a{\mathbf F}\;{\bigcirc_{s, k}}\;T_b{\mathbf F}=T_c{\mathbf F}
 \end{equation}
\end{definition}
 By virtue of Theorem \ref{symmconvo} it follows that the following theorem holds.
 \begin{theorem}\label{equivdef}
 A p.m. $\mathbf G\mbox{ is } \bigcirc_{s, k}-ID$, resp. stable if and only if
 $H_{s, k}({\mathbf G})$ is ID, resp. stable,  in the usual sense.
 \end{theorem}
  The following lemma will be used in the representation of $\bigcirc_{s, k}-ID, k\ge 2.$
    \begin{lemma}\label{Bessellimittheorem}
   (i) For every $t\ge 0$
   \begin{equation}\label{Bessellimittheorem 1}
   \lim_{x\rightarrow 0}\frac{1-\Lambda_s(tx)}{x^2}=
   \lim_{x\rightarrow
   0}\frac{1-Ee^{it\theta}}{x^2}=\frac{t^2}{2}.
   \end{equation}
 (ii) For any ${\mathbf x}=(x_0, x_1,\cdots ,x_k)\mbox{ and }{\mathbf t}=(t_0, t_1, \cdots, t_k)\in\mathbb R^{k+1}, k=1,2, ...$
 \begin{equation}\label{bessellimittheorem2}
 lim_{\rho\rightarrow 0}\frac{1-\prod_{r=0}^k \Lambda_s
(t_rx_r)}{\rho^2}=\frac{1}{2}\Sigma_{r=0}^k \lambda^2_r( Arg({\mathbf x}))t_r^2,
 \end{equation}
 where $\rho=||\mathbf x||, Arg({\mathbf x})=\frac{\mathbf x}{||\mathbf x||},\mbox{ and } \lambda_r( Arg({\mathbf x})), r=0,1, ...,k$ are given by
 \begin{equation}\label{polarization}
       \lambda_r( Arg({\mathbf x}))=
       \begin{cases}\cos\phi & r=0,\\
         \sin\phi\sin\phi_1\cdots
   \sin\phi_{r-1}\cos\phi_{r} &1\le r\le k-2,\\
   \sin\phi\sin\phi_1...\sin\phi_{k-2}\cos\psi & r={k-1},\\
   \sin\phi\sin\phi_1
   ...\sin\phi_{k-2}\sin\psi & r=k,
   \end{cases}
   \end{equation}
   where $0\le \psi, \phi, \phi_r\le\pi/2, r=1,2,...,k-2$ are angles
   of $\mathbf{x}$ appearing its polar form.
         \end{lemma}
 The following theorem gives a representation of rad.ch.f.'s of  $\bigcirc_{s, k}-$ID distributions
 (see \cite{Ng4} ), Theorem 2.6  for the proof).
 \begin{theorem}\label{LevyID} A p.m. $\mu\in ID(\bigcirc_{s, k})$ if and only if
there exist a $\sigma$-finite measure M (a L\'evy's measure) on
  $ \mathbb R^{+k}$ with the property that $M({\mathbf 0})=0,  {\mathbf M}$  is finite outside every neighborhood of ${\mathbf 0}$ and
\begin{equation}\label{integrable w. r. t. weight function}
\int_{\mathbb R^{+k}}\frac {\|{\mathbf x}\|^2} {1+\|{\mathbf x}\|^2}
{\mathbf M}(d{\mathbf x}) < \infty
\end{equation}
 and for each ${\mathbf t}=(t_1,...,t_k)\in
\mathbb R^{+k}$
\begin{equation}\label{Levy-Kintchine for k-dim.rad. ch. f.}
 -\log{\hat{\mu}({\mathbf  t})}=\int_{\mathbb R^{+k}}\{1-\prod_{j=1}^{k}\Lambda_s(t_jx_j)\} \frac
{1+\|{\mathbf x}\|^2} {\|{\mathbf x}\|^2} M({\mathbf {dx}}),
\end{equation}
where, at the origin $\mathbf{0}$, the integrand on the right-hand
side of (\ref{Levy-Kintchine for k-dim.rad. ch. f.}) is assumed to
be
\begin{equation}\label{limiting integrand}
\Sigma_{j=1}^k \lambda^2_j t_j^2 = lim_{\|\mathbf
{x}\|\rightarrow 0 }\{1-\prod_{j=1}^k
 \Lambda_s(t_jx_j)\} \frac {1+\|\mathbf x\|^2}
{\|\mathbf {x}\|^2}
\end{equation}
for nonnegative $\lambda_j, j=1, 2,...,k$  given by equations (\ref{polarization})
in Lemma \ref{Bessellimittheorem}.
 In particular, if $ M=0, \mbox{ then } \mu $ becomes a Rayleighian distribution with
the rad.ch.f (see definition \ref{Rayleigh})
\begin{equation}\label{kRayleighian rad. ch. f.}
-\log{\hat{\mu}({\mathbf t})}=\frac{1}{2}\sum_{j=1}^k \lambda^2_j
t_j^2,\quad {\mathbf t}\in \mathbb R^{+k},
\end{equation}
 for some nonnegative $\lambda_j, j=1,...,k.$
 Moreover, the representation (\ref{Levy-Kintchine for k-dim.rad. ch.
 f.}) is unique.
  \end{theorem}
   An immediate consequence of the above theorem is the following:
\begin{corollary}\label{Cor:Pair}
Each distribution $\mu\in ID(\bigcirc_{s, k})$ is uniquely determined by the pair $[\mathbf{M}, \pmb {\lambda}]$, where $\mathbf{M}$ is a  Levy's measure
on $\mathbb R^{+k}$ such that $\mathbf{M}(\mathbf{0})=0,$ $\mathbf{M}$ is finite outsite every neighbourhood of $\mathbf{0}$ and the condition (\ref{integrable w. r. t. weight function})
 is satisfied and $\pmb{\lambda}:=\{\lambda_1, \lambda_2,\cdots \lambda_k\}\in \mathbb R^{+k}$ is a vector of nonnegative numbers appearing in (\ref{kRayleighian rad. ch. f.}).
Consequently, one can write $\mu\equiv[\mathbf{M}, \pmb {\lambda}].$\\ \indent
In particular, if $\mathbf{M}$ is zero measure then $\mu=[\pmb{\lambda}]$ becomes a Rayleighian p.m. on $\mathbb R^{+k}$ as defined as follows:
\end{corollary}
  \begin{definition}\label{Rayleigh}
  A  k-dimensional distribution, say $\pmb{\mathbf \Sigma}_{s, k}$, is called  a {\it Rayleigh distribution}, if   \begin{equation}\label{k-dimension Rayleigh}
\pmb{\mathbf  \Sigma}_{s, k}=\sigma_s\times\sigma_s\times\cdots\times\sigma_s \quad
 (k\;times).
 \end{equation}
 Further, a distribution ${\mathbf F}\in \mathcal P(\mathbb R^{+k})$ is called a {\it Rayleighian distribution} if there exist nonnegative numbers $\lambda_r,
   r=1,2 \cdots k $ such that
\begin{equation}\label{k-dimensional rayleighian}
{ \mathbf F}=\{T_{\lambda_1}\sigma_s\}\times \{T_{\lambda_2}\sigma_s\}
 \times\ldots \times\{T_{\lambda_k}\sigma_s\}.
 \end{equation}
 \end{definition}
 \indent
 It is evident that every Rayleigh distribution is  a Rayleighian distribution. Moreover, every Rayleighian distribution is $\bigcirc_{s, k}-$ID. By virtue of (\ref{Maxwell density} )  and (\ref{k-dimension Rayleigh}) it follows that the k-dimensional Rayleigh density is given by
 \begin{equation}\label{density k-dimension Rayleigh}
 g({\mathbf x})=\Pi_{j=1}^k\frac{2^k(s+1)^{k(s+1)}}{\Gamma^k(s+1)}x_j^{2s+1}exp\{-(s+1)||{\mathbf x}||^2\},
 \end{equation}
 where ${\mathbf x}=(x_1, x_2,\ldots, x_k)\in \mathbb R^{+k}$ and the corresponding rad.ch.f. is given by
 \begin{equation}
 \hat\Sigma_{s, k}({\mathbf t})=Exp(-|{\mathbf t}|^2/4(s+1)),\quad {\mathbf t}\in \mathbb R^{+k}.
 \end{equation}
 Finally, the rad.ch.f. of a Rayleighian distribution $\mathbf  F\mbox{ on } \mathbb R^{+k}$ is given by
 \begin{equation}\label{rad.ch.rayleighian}
 \hat{\mathbf F}({\mathbf t})=Exp(-\frac{1}{2}\sum_{j=1}^k\lambda_j^2t_j^2)
 \end{equation}
 where $\lambda_j, j=1, 2, \ldots, k$ are some nonnegative numbers.
\section{An analogue of the L\'evy-Cram\'er Theorem in multi-dimensional Kingman convolution algebras}
  We say that a distribution ${\mathbf F \mbox{ on } \mathbb R^k}$ has dimension m, $1\le m \le k$,
 if m is the dimension of the smallest hyper-plane which contains the support of $\mathbf F.$
 The following theorem can be regarded as a version of the L\'evy-Cram\'er  Theorem
 for multi-dimensional Kingman convolution.The case k=1 was proved by Urbanik (\cite{U2}).
 \begin{theorem}\label{Levy-Cramer}
 Suppose that $\mathbf G_i \in \mathcal P(\mathbb R^{+k}),  i=1, 2 $  and
 \begin{equation}\label{decomposi}
\Sigma_{s, k}={\mathbf G}_1 \bigcirc_{s, k} {\mathbf G}_2.
 \end{equation}
 Then, ${\mathbf G}_i, i=1, 2$ are both Rayleighian distributions fufilling the condition that there exist nonnegative numbers $\lambda_{i, r}, i=1, 2\mbox{ and }  r=1, 2,\ldots, k$
 such that for each i=1, 2 the number of non-zero coefficients ${\lambda_{i, r}}'s$ among $\lambda_{i, 1}, \lambda_{i, 2},\ldots, \lambda_{i, k}$ are equal to the dimension of ${\mathbf G}_i,$ respectively.  Moreover,
 \begin{equation}
  \lambda_{1, r}^2+\gamma_{2, r}^2=1,\qquad r=1, 2, ..., k
  \end{equation}
 and
 \begin{equation}\label{form i}
{ \mathbf G}_i =T_{\lambda_{i, 1}}\sigma_s\times T_{\lambda_{i, 2}}\sigma_s\times \ldots\times T_{\lambda_{i, k}} \sigma_s
 \end{equation}

 \end{theorem}
 \begin{proof}
	 Suppose that  the equation (\ref{decomposi}) holds. Using the map $F_{s, k}$ we have
\begin{equation*}
F_{s, k}(\Sigma_{s, k})=F_{s, k}({\mathbf G}_1)\big.\ast F_{s, k}({\mathbf G}_2)
\end{equation*}	
which, by virtue of (\ref{Gaussian-Rayleigh}), implies that
\begin{equation*}
N({\mathbf 0}, 2(s+1){\mathbf I})=F_{s, k}({\mathbf G}_1)\big.\ast F_{s, k}({\mathbf G}_2).
\end{equation*}
  By the well-known L\'evy-Cram\'er Theorem on $\mathbb R^k$ (cf. Linnik and Ostrovskii \cite{LiOst}), that they are both symmetric Gaussian distributions on $\mathbb R^k.$ Consequently, they must be of the form
 (\ref{form i}) and the coefficients ${\lambda_{i, r}}'s$ satisfy
 the above stated conditions.
      \end{proof}
         
\end{document}